# Geometric and algebraic interpretation of primitive Pythagorean triples parameters

*Aleshkevich NataliaV.*

*Peter the Great St. Petersburg Polytechnic University*

**Abstract**

The paper found a geometric and algebraic interpretation of the parameters **m** and **n** from the formulas for obtaining primitive Pythagorean triples, which are solutions of the equation $x^2 + y^2 = z^2$:

$$x = m^2 - n^2, \ y = 2mn, \ z = m^2 + n^2.$$

The study was based on the process of building figurate numbers using gnomons. The paper discusses the process of building squares. The addition of the gnomon **U** to the original square leads to a larger square:

$$x^2 + U = z^2.$$

The first stage of the investigation was the construction of the gnomon **U**, which is equal to the area of a square **y²**. The construction is based on a generating square with a side equal to an even number. The area of the generating square is represented as the sum of the areas of two equal rectangles in all possible ways. At the same time, using the generating square, the gnomon **U** and the sides of all squares are also constructed: **x, y, z**.

The second stage of the investigation was to obtain a formula for the parameters **m** and **n** through the partition elements **t** and **l** of the side of the generating square.

## Introduction

The objects of our consideration are the Pythagorean numbers, also called primitive Pythagorean triples – triples **(x, y, z)** of positive integers satisfying the Pythagorean equation

$$x^2 + y^2 = z^2. \quad (1)$$

The Pythagorean Theorem is a fundamental geometric statement: in any right triangle the square of the hypotenuse is equal to the sum of the squares of the other two sides. This equation belongs to the class of Diophantine equation of second degree. The solution is sought in the form of integers.

The source of formulas that allow ancient Greek mathematicians to find different sets of Pythagorean triples, served as certain arithmetic identities [1, 158-161].

Thus, Pythagoras, using the identity

$$(2n + 1)^2 + (2n^2 + 2n)^2 = (2n^2 + 2n + 1)^2, \quad (2)$$

indicated formulas

$$(2n + 1, 2n(n + 1), 2n(n + 1) + 1), \quad (3)$$

describing all triples containing two consecutive numbers, one of which is the hypotenuse.

Plato, by using the identity

$$(m^2 - 1)^2 + (2m)^2 = (m^2 + 1)^2, \quad (4)$$

created the formulas

$$(m^2 - 1, 2m, m^2 + 1), \quad (5)$$

Which allowed to define, for even **m**, all primitive Pythagorean triples, in which hypotenuse and one of cathetus - consecutive odd numbers.

Euclid and Diophant, using the identity

$$(m^2 - n^2)^2 + (2mn)^2 = (m^2 + n^2)^2 \quad (6)$$

for relatively prime numbers **m>n** looked at following formulas

$$(m^2 - n^2, 2mn, m^2 + n^2), \quad (7)$$

described each Pythagorean triple *(x, y, z)* satisfying the condition $GCD(x, y) \leq 2$.

Although each of these formulas allows to get infinitely many Pythagorean triples, none of them is universal, that means it cannot create all possible triples.

We will consider triangles whose sides are expressed by relatively prime numbers, i.e., the numbers whose greatest common divisor is 1. Pythagorean triple *(x, y, z)* satisfying condition $GCD(x, y, z) = 1$ is called primitive.

In any primitive Pythagorean triple, one of the cathetus is an even number, and the other is an odd number. In this case, the hypotenuse *z* is an odd number. Without loss of generality, we assume that *x* is odd and *y* is even. With the above limitations, we can get all the primitive Pythagorean triples, and only them.

The best known in this area is the following statement.

*All primitive Pythagorean triples (x, y, z), for which y is an even number, can be obtained from the formulas:*

$$x = m^2 - n^2, y = 2mn, z = m^2 + n^2, m > n, (8)$$

*where m and n are relatively prime positive integers of different parity. Each primitive Pythagorean triple with even y is determined by this method uniquely.*

The parameters *m* and *n*, which form the primitive Pythagorean triples, were derived from rather abstract considerations. Until now, no geometric or algebraic interpretation of these quantities has been given.

There is also number of ways to describe Pythagorean triangles using flatness points.

One of them is based on the transition from equation $a^2 + b^2 = c^2, a, b, c \in N$ to equation $(a/c)^2 + (b/c)^2 = 1, a/c, b/c \in Q$ [2, 69-71]. This transition suggests that the set of primitive Pythagorean triangles is in one-to-one correspondence with the set of points with rational positive coordinates belonging to the unit circle $x^2 + y^2 = 1$. The general view of all rational solutions of the equation will be

$$x = \frac{m^2 - n^2}{m^2 + n^2} \quad (9)$$

$$y = \frac{2mn}{m^2 + n^2} \quad (10)$$

Introducing the notation $a = m^2 - n^2, b = 2mn, c = m^2 + n^2$, we obtain a primitive Pythagorean triple **(a, b, c)** corresponding to the rational point $(x, y), x = \frac{a}{c}, y = \frac{b}{c}$ of the unit circle. As we see, in this solution the parameters **m** and **n** are not physically defined.

Our task to find the process of obtaining the numbers **m** and **n** and discover their internal structure.

**Generation of parameters *m* and *n* for primitive Pythagorean triples**

The basis of our study will be the figurate numbers. Figurate numbers are mathematical objects that unite the concepts of number and geometric figure. The geometric figure appears to be a regular polygon. Each polygon is obtained from the previous one by adding the corresponding gnomon. Gnomon is an L-shaped figure, the application of which to the main structure does not change its shape. In our study, the geometric figure will be a square.

We will consider the construction the sequence of squares, starting with square has a side length equal to 1. The formula for this construction is:

$$(n + 1)^2 = n^2 + 2n + 1. \quad (11)$$

To each previous square we add the double value of its side plus 1.

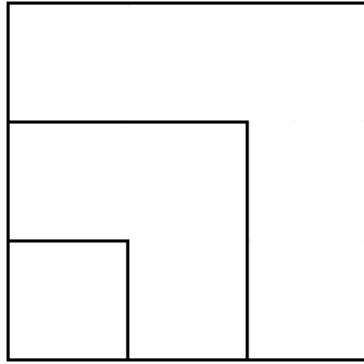
Picture 1

As a result, a "corner" is formed, which is the gnomon. Landing a gnomon on a square, or in other words, a square plus a gnomon will also be a square

$$x^2 + U = z^2. (12)$$

We need to build a gnomon so that it is equal in area to some square: $U = y^2$.

Then we come to the equation

$$x^2 + y^2 = z^2. (13)$$

## 1. Creation of the gnomon, who has the area of the square

Take a square with a side equal to an even number. Its side equals $S = 2tl$. The number **l** is odd number. The number **t** can be of any even number. In this case $GCD\,(t, l) = 1$.

The area of the square will be equal

$$S^2 = 2 \times 2t^2 \times l^2. (14)$$

Half of the area of the square will be equal to

$$\frac{S^2}{2} = 2t^2 \times l^2. (15)$$

Let's construct a square with side **y**

$$y = S + 2t^2 = 2tl + 2t^2 = 2t(l + t). (16)$$

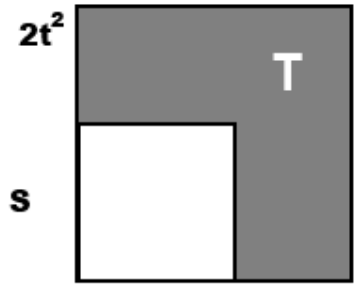

Picture 2

This side is the side of the gnomon **T** with a thickness equal $2t^2$. Let's lengthen the side of the gnomon **T** by the value of $l^2$.

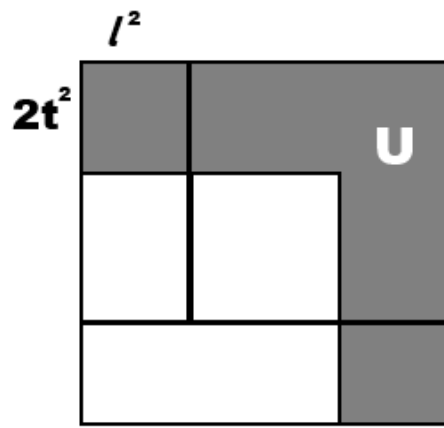

Picture 3

So, we built a gnomon **U** with a side

$$2t(l + t) + l^2 \quad (17)$$

and thickness $2t^2$.

The area of the gnomon **U** in this case will be equal to the area of the square with the side **y**.

We cut out the area of the original square with side **S**, from the area of the newly constructed square with side **y**, had previously converted this area into two identical rectangles with sides $2t^2$ and $l^2$ each. And placed these rectangles on two free sides of the gnomon **T**, thus having finished it to a gnomon **U**.

The side of the gnomon **U** is an odd number

$$z = 2tl + 2t^2 + l^2. \quad (18)$$

When subtracting the thickness of the gnomon from this side, we get the side of the second summable square

$$x = 2tl + l^2. \quad (19)$$

The area of the gnomon is equal to the area of the square with side **y**.

So, we have constructed three squares with the ratio

$$x^2 + y^2 = z^2. \quad (20)$$

Replacing **y**, **z** and **x** in accordance with formulas (16), (18), (19), we obtain the formula for the sum of the squares through parameters **t** and **l**:

$$(2tl + l^2)^2 + (2tl + 2t^2)^2 = (2tl + 2t^2 + l^2)^2. \quad (21)$$

## 2. The number of possible partitions of the side of the generating square

Each partition of the side of the generating square corresponds to one primitive Pythagorean triple.

The number 2 will be a separate constant. The remaining power of two $2^{\alpha_0 - 1}$ will be a multiplier of **t**. When $\alpha_0 = 1$, **t** is an odd number.

Prime odd factors for **t** are chosen from the product. The remaining unselected factors will be **l**. Each prime factor is taken with its power. Wherein $GCD(t, l) = 1$.

You can select zero odd factors. It corresponds $C_r^0$. You can select any single item. The number of samples is $C_r^1$. You can choose any two items. The number of samples is $C_r^2$. And so on to **r** elements. In the latter case, the number of samples will be $C_r^r$.

The total number of different samples from **r** is $\sum_{i=0}^{r} C_r^i = 2^r$.

Whatever number of odd factors is chosen from the factors of the number **S** for the parameter **t**, the remaining factors will be included in the parameter **l**.

### 3. Expression of the parameters *m* and *n* through the partitioning parameters of the side of the generating square

Let's make the following transformation:

represent **y** as $2t(l + t)$.

We choose that **l** is an odd number. It can be either 1 or the product of odd prime factors in the corresponding powers. The number **t** can be any parity. It can be either 1 or the product of prime factors in the corresponding powers. We get that the factors **t** and **(l+t)** have different parity. And, indeed, if **t** is even, then the contents of the bracket are odd and vice versa.

The transition from our notation to the generally accepted in terms of **m** and **n**:

$$m = l + t; n = t. \quad (22)$$

In this case, one of the summable squares will have a side $2mn$. The total square will have a side $m^2 + n^2$. In fact, $2tl + 2t^2 + l^2 = t^2 + 2tl + l^2 + t^2 = (t + l)^2 + t^2$.

The second square will have a side $m^2 - n^2$. In fact, $2tl + l^2 = 2tl + l^2 + t^2 - t^2 = (t + l)^2 - t^2$.

The basis of our construction is a square with an even side, generating the sum of three squares. A sequence of squares, starting with a square with a side equal to 2 (with the rules stated above, for partition the side of a square into two subsets, one of which has an even value and the second one is odd), leads to the construction of all possible primitive Pythagorean triples.

As a result of our investigation were revealed the geometrical interpretation and algebraic meaning that underlie the choice of the parameters **m** and **n**.

**m** - is the sum of the values of two subsets of the partition;

**n** - is the value of one subset of the partition, it necessarily includes the factor $2^{\alpha_0-1}$ with $\alpha_0 > 1$.

The number of **L(S)** partitions

$S = 2 \times 2^{\alpha_0-1} p_1^{\alpha_1} p_2^{\alpha_2} \ldots p_r^{\alpha_r}$ equals $L(S) = 2^r$. (23)

**Conclusion**

We have found a geometric and algebraic interpretation for the parameters **m** and **n** from the general formulas for the definition of primitive Pythagorean triples. This interpretation is based on the concept of a figurate number and the way figurate number is constructed with the help of a gnomon.

For any even number, starting with 2, it always possible to construct a primitive Pythagorean triple. The number of such constructions depends on the representation of the number as a partition of the product of prime factors with the corresponding powers. The number of such constructions depends only on the number of odd prime factors without taking their degrees into account. The number of **L(S)** of partition of the product

$S = 2 \times 2^{\alpha_0-1} p_1^{\alpha_1} p_2^{\alpha_2} \ldots p_r^{\alpha_r}$ equals $L(S) = 2^r$.

**m** - is the sum of two elements of the partitioning of the product of an even number (**t, l**):

$$m = t + l;$$

**n** - is one of the elements of the partition, namely:

$$n = t.$$

**l** - is always odd.